\documentstyle[12pt, amsfonts]{amsart}
\begin{document}
\def\R{{\mathbb R}}
\def\Z{{\mathbb Z}}
\def\C{{\mathbb C}}
\newcommand{\trace}{\rm trace}
\newcommand{\Ex}{{\mathbb{E}}}
\newcommand{\Prob}{{\mathbb{P}}}
\newcommand{\E}{{\cal E}}
\newcommand{\F}{{\cal F}}
\newtheorem{df}{Definition}
\newtheorem{theorem}{Theorem}
\newtheorem{lemma}{Lemma}
\newtheorem{pr}{Proposition}
\newtheorem{co}{Corollary}
\newtheorem{problem}{Problem}
\def\n{\nu}
\def\sign{\mbox{ sign }}
\def\a{\alpha}
\def\N{{\mathbb N}}
\def\A{{\cal A}}
\def\L{{\cal L}}
\def\X{{\cal X}}
\def\F{{\cal F}}
\def\c{\bar{c}}
\def\v{\nu}
\def\d{\delta}
\def\diam{\mbox{\rm dim}}
\def\vol{\mbox{\rm Vol}}
\def\b{\beta}
\def\t{\theta}
\def\l{\lambda}
\def\e{\varepsilon}
\def\colon{{:}\;}
\def\pf{\noindent {\bf Proof :  \  }}
\def\endpf{ \begin{flushright}
$ \Box $ \\
\end{flushright}}

\title[Slicing inequalities for measures]
{Slicing inequalities for measures of convex bodies}

\author{Alexander Koldobsky}

\address{Department of Mathematics\\ 
University of Missouri\\
Columbia, MO 65211}

\email{koldobskiya@@missouri.edu}

\begin{abstract}  
We consider the following problem. Does there exist an absolute constant $C$ so that for every $n\in \N,$
every integer $1\le k < n,$ every origin-symmetric convex body $L$ in $\R^n,$ and every measure $\mu$ with non-negative even continuous density in $\R^n,$
\begin{equation}\label{mainproblem}
\mu(L)\ \le\ C^k
 \max_{H \in Gr_{n-k}} \mu(L\cap H)\  |L|^{k/n},
\end{equation}
where  $Gr_{n-k}$ is the Grassmanian of $(n-k)$-dimensional subspaces of $\R^n$, 
and $|L|$ stands for volume? This question is an extension to arbitrary measures (in place of volume) and
to sections of arbitrary codimension $k$ of the slicing problem of Bourgain, a major open problem 
in convex geometry. 

It was proved in \cite{K4,K5} that (\ref{mainproblem}) holds for arbitrary origin-symmetric convex bodies, 
all $k$ and all $\mu$ with $C\le O(\sqrt{n}).$ In this article, we prove inequality (\ref{mainproblem}) with an absolute 
constant $C$ for unconditional convex bodies and for duals of bodies with bounded volume ratio. We also prove
that for every $\lambda\in (0,1)$ there exists a constant $C=C(\lambda)$  so that inequality (\ref{mainproblem}) 
holds for every $n\in \N,$ every origin-symmetric convex body $L$ in $\R^n,$ every measure $\mu$ with continuous density
and the codimension of sections $k\ge \lambda n.$ 
The proofs are based on a stability result for generalized intersection bodies
and on estimates of the outer volume ratio distance from an arbitrary convex body to the classes
of generalized intersection bodies. In the last section, we show that for some measures the behavior 
of minimal sections may be very different from the case of volume. 
\end{abstract}  
\maketitle

\section{Introduction}
The slicing problem \cite{Bo1, Bo2, Ba1, MP}, a major open problem in convex geometry,
asks whether there exists an absolute constant $C$ so that for any origin-symmetric convex body $K$ in $\R^n$
of volume 1 there is a hyperplane section of $K$ whose $(n-1)$-dimensional volume is greater than $1/C.$
In other words, does there exist a constant $C$ so that for any $n\in \N$ and any
origin-symmetric convex body $K$ in $\R^n$
\begin{equation} \label{hyper}
|K|^{\frac {n-1}n} \le C \max_{\xi \in S^{n-1}} |K\cap \xi^\bot|,
\end{equation}
where  $\xi^\bot$ is the central hyperplane in $\R^n$ perpendicular to $\xi,$ and
$|K|$ stands for volume of proper dimension?
The best current result $C\le O(n^{1/4})$ is due to Klartag \cite{Kl}, who
removed the  logarithmic term from an earlier estimate of Bourgain \cite{Bo3}.
We refer the reader to [BGVV] for the history and partial results.

For certain classes of bodies the question has been answered in affirmative. These classes
include unconditional convex bodies (as initially observed by Bourgain; see also \cite{MP, J2,
BN, BGVV}), unit balls of subspaces of $L_p$  \cite{Ba2, J1, M1}, intersection bodies
\cite[Theorem 9.4.11]{G}, zonoids, duals of bodies with bounded volume ratio
\cite{MP}, the Schatten classes \cite{KMP}, $k$-intersection bodies \cite{KPY, K6}.

Iterating (\ref{hyper}) one gets the lower dimensional slicing problem asking whether
the inequality
\begin{equation} \label{lowdimhyper}
|K|^{\frac {n-k}n} \le C^k \max_{H\in Gr_{n-k}} |K\cap H|
\end{equation}
holds with an absolute constant $C$ where $1\le k \le n-1$ and $Gr_{n-k}$
is the Grassmanian of $(n-k)$-dimensional subspaces of $\R^n.$

In this note we prove (\ref{lowdimhyper}) in the case where $k\ge \lambda n,\ 0<\lambda<1,$
with the constant $C=C(\lambda)$ dependent only on $\lambda.$ Moreover, we prove this result in
a more general setting of arbitrary measures in place of volume. We
consider the following generalization of the slicing problem.

\begin{problem} \label{prob}
Does there exist an absolute constant $C$ so that for every $n\in \N,$
every integer $1\le k < n,$ every origin-symmetric convex body $L$ in $\R^n,$ and every measure $\mu$ with non-negative even continuous density $f$ in $\R^n,$
\begin{equation}\label{main-problem}
\mu(L)\ \le\ C^k  \max_{H \in Gr_{n-k}} \mu(L\cap H)\ |L|^{k/n}.
\end{equation}
\end{problem}
Here $\mu(B)=\int_B f$ for every compact set $B$ in $\R^n,$ and $\mu(B\cap H)=\int_{B\cap H} f$
is the result of integration of the restriction of $f$ to $H$ with respect to Lebesgue measure in $H.$

In many cases we will write (\ref{main-problem}) in an equivalent form
\begin{equation}\label{measslicing}
\mu(L)\ \le\ C^k \frac n{n-k}\ c_{n,k} \max_{H \in Gr_{n-k}} \mu(L\cap H)\ |L|^{k/n},
\end{equation}
where $c_{n,k}= |B_2|^{\frac {n-k}n}/|B_2^{n-k}|,$ and $B_2^n$ is the unit Euclidean ball in $\R^n.$
Note that $c_{n,k}\in (e^{-k/2},1)$ (see for example \cite[Lemma 2.1]{KL}), and
$$1\le  \frac{n}{n-k} \le e^{\frac k{n-k}}\le  e^k,$$ so these constants
can be incorporated in the constant $C.$ 

It appears that some results on the original slicing problem can be extended to the case of arbitrary measures.
The first result of this kind was established in \cite{K3}, namely, when $L$ is an intersection body (see
definition below) and $k=1,$ inequality (\ref{measslicing}) holds with the best possible constant $C=1.$
This result was later proved for arbitrary $k$ in \cite{KM}.
For arbitrary origin-symmetric convex bodies, inequality (\ref{measslicing}) was proved with $C=\sqrt{n}$ 
in \cite{K4} and \cite{K5}, for $k=1$ and for arbitrary $k$, respectively.
When $L$ is the unit ball of a subspace of $L_p,\ p\ge 2,$ the constant $C$ can be improved to
$n^{\frac 12-\frac 1p};$ see \cite{K6}. 
In \cite{K6}, (\ref{main-problem}) was also proved for the unit balls of normed spaces that embed
in $L_p,\ -\infty<p\le 2$ with $C$ depending only on $p.$ In the case where $k=1$ and the measure $\mu$
is log-concave, (\ref{main-problem}) holds for any origin-symmetric convex body with
$C\le O(n^{1/4}),$ as shown in \cite{KZ} using the estimate of Klartag \cite{Kl} mentioned above
and the technique of Ball \cite{Ba1} relating log-concave measures to convex bodies.

In this article, we prove inequality (\ref{main-problem}) for unconditional convex bodies and duals
of bodies with finite volume ratio, with an absolute constant $C.$ We also prove
that for every $\lambda\in (0,1)$ there exists a constant $C=C(\lambda)$  so that 
inequality (\ref{main-problem}) holds for every $n\in \N,$ arbitrary origin-symmetric convex 
body $L,$ every measure $\mu$ with continuous density
and every codimension of sections $k$ satisfying $\lambda n\le k <n.$ 

In Section \ref{counter}, we show that the properties of the minimal measures of sections
may be different from the case of volume. We prove that there exist a symmetric convex
body $L$ in $\R^n$ and a measure $\mu$ with continuous density so that
$$\mu(L)< \frac n{n-1} c_{n,1} \min_{\xi\in S^{n-1}} \mu(L\cap \xi^\bot) |L|^{1/n}.$$
Note that in the case of volume 
$$ \int_{S^{n-1}} |K\cap \xi^\bot| d\sigma(\xi) \le c_{n,1} |K|^{\frac{n-1}n},$$
where $\sigma$ is the normalized uniform measure on the sphere; see \cite{L1} for more 
general results.
 
\section{Reduction to intersection bodies} \label{slicing}

The approach to Problem \ref{prob} suggested in this paper is based on the concept
of an intersection body. In this section we reduce the problem to computing the
outer volume ratio distance from an origin-symmetric convex body to the class
of generalized intersection bodies.

We need several definitions and facts.
A closed bounded set $K$ in $\R^n$ is called a {\it star body}  if 
every straight line passing through the origin crosses the boundary of $K$ 
at exactly two points different from the origin, the origin is an interior point of $K,$
and the {\it Minkowski functional} 
of $K$ defined by 
$$\|x\|_K = \min\{a\ge 0:\ x\in aK\}$$
is a continuous function on $\R^n.$ 

The {\it radial function} of a star body $K$ is defined by
$$\rho_K(x) = \|x\|_K^{-1}, \qquad x\in \R^n,\ x\neq 0.$$
If $x\in S^{n-1}$ then $\rho_K(x)$ is the radius of $K$ in the
direction of $x.$

We use the polar formula for volume of a star body
\begin{equation}\label{polar}
|K|=\frac 1n \int_{S^{n-1}} \|\theta\|_K^{-n} d\theta.
\end{equation}

The class of intersection bodies was introduced by Lutwak \cite{L2}.
Let $K, L$ be origin-symmetric star bodies in $\R^n.$ We say that $K$ is the 
intersection body of $L$ if the radius of $K$ in every direction is 
equal to the $(n-1)$-dimensional volume of the section of $L$ by the central
hyperplane orthogonal to this direction, i.e. for every $\xi\in S^{n-1},$
$$
\rho_K(\xi)= \|\xi\|_K^{-1} = |L\cap \xi^\bot|
$$
$$= \frac 1{n-1} \int_{S^{n-1}\cap \xi^\bot} \|\theta\|_L^{-n+1}d\theta=
\frac 1{n-1} R\left(\|\cdot\|_L^{-n+1}\right)(\xi),$$
where $R:C(S^{n-1})\to C(S^{n-1})$ is the {\it spherical Radon transform}
$$Rf(\xi)=\int_{S^{n-1}\cap \xi^\bot} f(x) dx,\qquad \forall f\in C(S^{n-1}).$$
All bodies $K$ that appear as intersection bodies of different star bodies
form {\it the class of intersection bodies of star bodies}. A more general class of {\it intersection bodies} 
is defined as follows. If $\mu$ is a finite Borel measure on $S^{n-1},$ then the spherical Radon transform
$R\mu$ of $\mu$ is defined as a functional on $C(S^{n-1})$ acting by
$$(R\mu, f)=(\mu, Rf)=\int_{S^{n-1}} Rf(x) d\mu(x),\qquad \forall f\in C(S^{n-1}).$$
A star body $K$ in $\R^n$ is called an {\it intersection body} if $\|\cdot\|_K^{-1}=R\mu$
for some measure $\mu,$ as functionals on $C(S^{n-1}),$  i.e.
$$\int_{S^{n-1}} \|x\|_K^{-1} f(x) dx = \int_{S^{n-1}} Rf(x)d\mu(x),\qquad \forall f\in C(S^{n-1}).$$
Intersection bodies played a crucial role in the solution of the
Busemann-Petty problem and its generalizations; see \cite[Chapter 5]{K1}.

A generalization of the concept of an intersection body
was introduced by Zhang \cite{Z} 
in connection with the lower dimensional Busemann-Petty problem.
For $1\le k \le n-1,$  the {\it $(n-k)$-dimensional spherical Radon transform} 
$R_{n-k}:C(S^{n-1})\to C(Gr_{n-k})$  
is a linear operator defined by
$$R_{n-k}g (H)=\int_{S^{n-1}\cap H} g(x)\ dx,\quad \forall  H\in Gr_{n-k}$$
for every function $g\in C(S^{n-1}).$
\smallbreak
We say that
an origin symmetric star body $K$ in $\R^n$ is a {\it generalized $k$-intersection body}, 
and write $K\in {\cal{BP}}_k^n,$  if there exists a finite Borel non-negative measure $\mu$
on $Gr_{n-k}$ so that for every $g\in C(S^{n-1})$
\begin{equation}\label{genint}
\int_{S^{n-1}} \|x\|_K^{-k} g(x)\ dx=\int_{Gr_{n-k}} R_{n-k}g(H)\ d\mu(H).
\end{equation}
When $k=1$ we get the class of intersection bodies.
It was proved by Grinberg and Zhang \cite[Lemma 6.1]{GZ} that every intersection body in $\R^n$
is a generalized $k$-intersection
body for every $k<n.$ More generally, as proved later by E.Milman \cite{M2}, if  $m$ divides $k$, then every
generalized $m$-intersection body is a generalized $k$-intersection body.  Note that in \cite{Z,GZ}
generalized $k$-intersection bodies are called ``$i$-intersection bodies".

We need a stability result for generalized $k$-intersection bodies proved in \cite[Theorem 1]{K5}.
Here we present a slightly simpler version.  
\begin{theorem}\label{stab2}
Suppose that $1\le k \le n-1,$ $K$ is a generalized $k$-intersection body in $\R^n,$  $f$
is an even continuous non-negative function on $K,$ and $\e>0.$ If
$$
\int_{K\cap H} f \ \le \e,\qquad \forall H\in Gr_{n-k},
$$
then
$$
\int_K f\ \le \frac {n}{n-k}\ c_{n,k}\ |K|^{k/n}\e.
$$
Recall that $c_{n,k}\in (e^{-k/2},1).$
\end{theorem}
\pf Writing integrals in spherical coordinates we get
$$
\int_K f = \int\limits_{S^{n-1}}\left(\int\limits_0^{\|\theta\|^{-1}_K} r^{n-1} f(r\theta)\ dr\right) d\theta,
$$
and
$$\int_{K\cap H} f =  
\int_{S^{n-1}\cap H} \left(\int_0^{\|\theta\|_K^{-1}} r^{n-k-1}f(r\theta)\ dr \right)d\theta$$
$$
=R_{n-k}\left(\int_0^{\|\cdot\|_K^{-1}} r^{n-k-1}f(r\ \cdot)\ dr \right)(H),
$$
so the condition of the theorem can be written as
$$R_{n-k}\left(\int_0^{\|\cdot\|_K^{-1}} r^{n-k-1}f(r\ \cdot)\ dr \right)(H) \le  \e,\qquad \forall H\in Gr_{n-k}.$$
Integrate both sides with respect to the measure $\mu$ on $Gr_{n-k}$ that corresponds to $K$
as a generalized $k$-intersection body by (\ref{genint}). 
We get
$$\int_{S^{n-1}} \|\theta\|_K^{-k} \left(\int_0^{\|\theta\|_K^{-1}} r^{n-k-1}f(r\theta)\ dr \right)d\theta 
\le  \e \mu(Gr_{n-k}).$$
Estimate the integral in the left-hand side from below using $f\ge 0:$ 
$$\int_{S^{n-1}} \|\theta\|_K^{-k} \left(\int_0^{\|\theta\|_K^{-1}} r^{n-k-1}f(r\theta)\ dr \right)d\theta$$
$$=\int_{S^{n-1}}  \left(\int_0^{\|\theta\|_K^{-1}} r^{n-1}f(r\theta)\ dr \right)d\theta$$ 
$$+ \int_{S^{n-1}} \left(\int_0^{\|\theta\|_K^{-1}} (\|\theta\|_K^{-k} - r^k)  r^{n-k-1}f(r\theta)\ dr \right)d\theta$$
$$\ge \int_{S^{n-1}}  \left(\int_0^{\|\theta\|_K^{-1}} r^{n-1}f(r\theta)\ dr \right)d\theta= \int_Kf.$$
Now we estimate $\mu(Gr_{n-k})$ from above. We use $1= R_{n-k}1(H)/|S^{n-k-1}|$ for every $H\in Gr_{n-k},$ 
definition (\ref{genint}), H\"older's inequality and the fact that $n|B_2^n|=|S^{n-1}|$:
$$ \mu(Gr_{n-k})= \frac 1{\left|S^{n-k-1}\right|} \int_{Gr_{n-k}} R_{n-k}1(H) d\mu(H)$$
$$=\frac 1{\left| S^{n-k-1} \right| } \int_{S^{n-1}} \|\theta\|_K^{-k}\ d\theta $$
$$\le  \frac 1{\left|S^{n-k-1}\right|} \left|S^{n-1}\right|^{\frac{n-k}n} \left(\int_{S^{n-1}} \|\theta\|_K^{-n}\ d\theta\right)^{\frac kn}$$
$$
=  \frac{1}{\left|S^{n-k-1}\right|} \left|S^{n-1}\right|^{\frac{n-k}n} n^{k/n}|K|^{k/n}= \frac n{n-k} c_{n,k} |K|^{k/n}.
$$
Combining the estimates,
$$\int_K f\ \le\ \frac n{n-k} c_{n,k} |K|^{k/n} \e.$$
\endpf
\bigbreak
For a convex body $L$ in $\R^n$ and $1\le k <n,$ denote by 
$${\rm {o.v.r.}}(L,{\cal{BP}}_k^n) = \inf \left\{ \left( \frac {|K|}{|L|}\right)^{1/n}:\ L\subset K,\ K\in {\cal{BP}}_k^n \right\}$$
the outer volume ratio distance from a body $L$ to the class ${\cal{BP}}_k^n.$

\begin{co} \label{lowdim} Let $L$ be an origin-symmetric star body in $\R^n.$ Then for any measure $\mu$
with even continuous density on $L$ we have
$$\mu(L)\le \left({\rm{ o.v.r.}}(L,{{\cal{BP}}_k^n})\right)^k\ \frac n{n-k} 
c_{n,k} \max_{H\in Gr_{n-k}} \mu(L\cap H)\ |L|^{k/n}.$$
\end{co}

\pf  Let $C>{\rm o.v.r.}(L,{\cal{BP}}_k^n),$ then there exists a body $K$ in ${\cal{BP}}_k^n$ such that $L\subset K$
and $|K|^{1/n}\le C\ |L|^{1/n}.$

Let $g$ be the density of the measure $\mu,$ and define a function $f$ on $K$ by
$f= g \chi_L,$ where $\chi_L$ is the indicator function of $L.$ Clearly, $f\ge 0$ everywhere on $K.$ Put 
$$\e=\max_{H\in Gr_{n-k}} \int_{K\cap H} f = \max_{H\in Gr_{n-k}} \int_{L\cap H} g=\max_{H\in Gr_{n-k}}\mu(L\cap H),$$
and apply Theorem \ref{stab2} to $f,K,\e$ ($f$ is not continuous, but we can do an easy approximation). 
We have
$$\mu(L)= \int_L g = \int_K f 
 \le \frac n{n-k} c_{n,k} |K|^{k/n}\max_{H\in Gr_{n-k}} \mu(L\cap H)$$
$$ \le C^k\ \frac n{n-k} c_{n,k} |L|^{k/n}\max_{H\in Gr_{n-k}} \mu(L\cap H).$$
The result follows by sending $C$ to ${\rm o.v.r.}(L,{\cal{BP}}_k^n).$
\endpf
\section{Unconditional bodies}

 Let $e_i,\ 1\le i\le n,$ be the standard basis of $\R^n.$ A star body $K$ in $\R^n$ is called unconditional 
if for every choice of real numbers $x_i$ and $\delta_i=\pm 1,\ 1\le i \le n$
we have 
$$\|\sum_{i=1}^n \delta_ix_i e_i \|_K = \|\sum_{i=1}^n x_i e_i \|_K.$$

\begin{theorem}\label{uncond} For every $n\in \N,$
every $1\le k<n,$ every unconditional convex body $L$ in $\R^n$
and every measure $\mu$ with even continuous non-negative density on $L$
\begin{equation} \label{uncond12}
\mu(L)\ \le\ e^k \frac{n}{n-k}c_{n,k}
 \max_{H \in Gr_{n-k}} \mu(L\cap H)\ |L|^{k/n}.
\end{equation}
\end{theorem}

\pf  By a result of Lozanovskii \cite{Lo} (see the proof in \cite[Corollary 3.4]{P}), there exists a linear
operator $T: \R^n\to \R^n$ so that 
$$T(B_\infty^n) \subset L \subset n T(B_1^n),$$ 
where $B_1^n$ and $B_\infty^n$ are the unit balls of the spaces $\ell_1^n$ and $\ell_\infty^n,$
respectively.
Let $K=nT(B_1^n).$ By \cite[Theorem 3]{K2} and the fact that a linear transformation of an intersection body
is an intersection body (see \cite{L2} or \cite[Theorem 1]{K2}), the body $K$ is an intersection body in $\R^n.$
By a result of Grinberg and Zhang \cite[Lemma 6.1]{GZ}, $K$ is a generalized $k$-intersection body
for every $1\le k<n.$

Since $|B_1^n| = 2^n/n!$ (see for example \cite[Lemma 2.19]{K1}), we have
$|K|^{1/n}\le 2e |\det T|^{1/n}.$ On the other hand, $|T(B_\infty^n)| = 2^n |\det T|,$ and $T(B_\infty^n)\subset L,$
so $|K|^{1/n} \le e\  |L|^{1/n}.$ Therefore, ${\rm o.v.r}(L,{\cal{BP}}_k^n)\le e.$ Now (\ref{uncond12}) follows from
Corollary \ref{lowdim}. \endpf

\section{Duals of bodies with bounded volume ratio}

The volume ratio of a convex body $K$ in $\R^n$ is defined by
$${\rm{v.r.}}(K) =\inf_E \left\{ \left(\frac{|K|}{|E|}\right)^{1/n}:\ E\subset K,\ E-{\rm ellipsoid}\right\}.$$
The following argument is standard and first appeared in \cite{BM} and \cite{MP}.
Let $K^\circ$ and $E^\circ$ be polar bodies of $K$ and $E,$ respectively. If $E$ is an ellipsoid,
then 
$$|E| |E^\circ|=|B_2^n|^2.$$
By the reverse Santalo inequality of Bourgain and Milman \cite{BM}, there exists an absolute constant
$c>0$ such that
$$\left(|K||K^\circ|\right)^{1/n} \ge \frac cn.$$
Combining these and using the asymptotics of $B_2^n$ we get that there exists an absolute
constant $C$ such that
$$ \left(\frac{|E^\circ|}{|K^\circ|}\right)^{1/n} \le C \left(\frac{|K|}{|E|}\right)^{1/n}.$$

\begin{theorem}\label{boundedvr} There exists an absolute constant $C$ such that for every $n\in \N,$
every $1\le k<n,$ every origin-symmetric convex body $L$ in $\R^n$
and every measure $\mu$ with even continuous non-negative density on $L$
$$
\mu(L)\ \le\ (C\ {\rm v.r}(L^\circ))^k \frac{n}{n-k}c_{n,k}
\max_{H \in Gr_{n-k}} \mu(L\cap H)\ |L|^{k/n}.
$$
\end{theorem}

\pf If $E$ is an ellipsoid, $E\subset L^\circ$, then the ellipsoid $E^\circ$ contains $L.$
Also every ellipsoid  is an intersection body as a linear image of the Euclidean ball,
so it is also a generalized $k$-intersection body for every $k.$
By the argument before the statement of the theorem, 
$${\rm o.v.r}(L,{\cal{BP}}_k^n)\le C\ {\rm v.r.}(L^\circ).$$
The result follows from Corollary \ref{lowdim}.  \qed
\bigbreak

\section{Sections of proportional dimensions}

The outer volume ratio distance from a general convex body to the class of
generalized $k$-intersection bodies was estimated in \cite{KPZ}.

\begin{pr} \label{kpz} (\cite[Theorem 1.1]{KPZ}) Let $L$ be an origin-symmetric convex 
body in $\R^n,$ and let $1\le k \le n-1.$ Then
$${\rm o.v.r.}(L,{\cal{BP}}_k^n) \le C_0 \sqrt{\frac{n}{k}}\left(\log\left(\frac{en}{k}\right)\right)^{3/2},$$
where $C_0$ is an absolute constant.
\end{pr}

\noindent {\bf Remark.} In \cite[Theorem 1.1]{KPZ}, the result was formulated with the logarithmic
term raised to the power 1/2 instead of 3/2. This happened because the proof in \cite[p.2705]{KPZ}
uses Corollary 3.2 which holds for $\alpha=1.$ However, the constant $\alpha$ used in the proof is
$\alpha= 2- \frac 1{\log(en/k)},$ so Corollary 3.2 should have been formulated for
this different value of $\alpha.$ We now correct this at the expense of an extra logarithmic term.

We use a result of Pisier \cite[Corollary 7.16]{P}, generalizing V.Milman's $M$-position. For two symmetric
convex bodies K and L in $\R^n,$ the covering number of $K$ by $L,$ denoted by $N(K,L),$ is
defined as the minimal number of translates of $L,$ with their centers in K, needed to cover K.

\begin{theorem}\label{pisier} (\cite[p.120]{P}) For every $\alpha\in (0,2)$ and every origin-symmetric convex
body $K$ in $\R^n,$ there exists a linear image $K_\alpha$ of $K$ such that
$$\max \{ N(K_\alpha,tB_2^n), N(B_2^n,tK_\alpha)\} \le \exp \left(\frac{cn}{t^\alpha(2-\alpha)}\right),$$
for every $t\ge 1,$ where $c$ is an absolute constant.
\end{theorem}
The constant $c/(2-\alpha)$ is not written precisely in Corollary 7.16 of \cite{P}, but it
can be established by combining Corollary 7.15 and the proofs of Theorems 7.13  and 7.11 in
the same book.
\smallbreak
Theorem \ref{pisier} implies a generalization of V.Milman's reverse Brunn-Minkowski inequality;
one can find this in \cite{P} as a combination of several results. We present a proof for the sake 
of completeness.

\begin{co}\label{repl} Let $\alpha\in [1,2),$ let $K$ be an origin-symmetric convex body in $\R^n,$
and let $K_\alpha$ be the position of $K$ established in Theorem \ref{pisier}. Then for every $t\ge 1,$
$$|K_\alpha+tB_2^n|^{1/n} \le 2e^c\ t|K_\alpha|^{1/n} 
\frac 1{2-\alpha} \exp \left(\frac{c}{t^\alpha(2-\alpha)}\right),$$
where $c$ is the same absolute constant as in Theorem \ref{pisier}.
\end{co}
\pf We first use the part of Theorem \ref{pisier} estimating $N(B_2^n,tK_\alpha).$
Put $t=(2-\alpha)^{-1/\alpha}$ in Theorem \ref{pisier}. Then
$$|B_2^n|^{1/n} \le t |K_\alpha|^{1/n} \left(N(B_2^n,tK_\alpha)\right)^{1/n}$$$$ \le
 (2-\alpha)^{-1/\alpha} e^c |K_\alpha|^{1/n}\le \frac{e^c}{2-\alpha} |K_\alpha|^{1/n}.$$
Now for every $t\ge 1$ we use the estimate for $N(K_\alpha, tB_2^n)$ from Theorem \ref{pisier}.
We have
$$\frac{|K_\alpha+tB_2^n|^{1/n}}{2t|K_\alpha|^{1/n}}\le  \frac {e^c}{2-\alpha}
\frac{|K_\alpha+tB_2^n|^{1/n}}{2t|B_2^n|^{1/n}}$$
$$\le \frac {e^c}{2-\alpha} \left(N(K_\alpha+tB_2^n,2tB_2^n)\right)^{1/n}$$$$\le
 \frac {e^c}{2-\alpha} \left(N(K_\alpha,tB_2^n)\right)^{1/n}
\le \frac {e^c}{2-\alpha}\exp \left(\frac{c}{t^\alpha(2-\alpha)}\right). \qed$$
In the proof of Theorem 1.1. in \cite[p.2705]{KPZ}, we have
$\alpha= 2-\frac 1{\log e\frac nk}$ and $t^\alpha(2-\alpha)=\frac {n}k,$ so 
$t\sim \sqrt{\frac nk \log(\frac {en}k)}.$ Then Corollary \ref{repl} implies
$$|K_\alpha + tB_2^n|^{1/n}\le c' \sqrt{\frac nk}\left(\log\left(\frac {en}k\right)\right)^{3/2} |K_\alpha|^{1/n},$$
where $c'$ is an absolute constant. Using this estimate in place of Corollary 3.2 in \cite[p.2705]{KPZ},
we get Proposition \ref{kpz}.
\medbreak
Proposition \ref{kpz} in conjunction with Corollary \ref{lowdim} implies
the following slicing inequality.
 
\begin{theorem} There exists an absolute constant $C_0$ such that for every $n\in \N,$
every $1\le k < n,$ every origin-symmetric convex body $L$ in $\R^n$
and every measure $\mu$ with even continuous non-negative density on $L$
$$
\mu(L)\ \le\  C_0^k  \left(\sqrt{\frac{n}{k}}\left(\log\left(\frac{en}{k}\right)\right)^{3/2}\right)^k \frac n{n-k}c_{n,k}
\max_{H\in Gr_{n-k}}  \mu(L\cap H)\ |L|^{k/n}.
$$
\end{theorem}

\begin{co} If the codimension of sections $k$ satisfies $\lambda n\le k <n,$ for some $\lambda\in (0,1),$ 
then for every origin-symmetric convex body $L$ in $\R^n$ and every measure $\mu$ with continuous
non-negative density in $\R^n,$
$$ \mu(L)\ \le\  C_0^k 
\left(\sqrt{\frac{(1-\log \lambda)^3}{\lambda}}\right)^k \frac{n}{n-k}c_{n,k}
 \max_{H \in Gr_{n-k}} \mu(L\cap H)\ |L|^{k/n},$$
where $C_0$ is an absolute constant.
\end{co}

\section{Minimal sections} \label{counter}

We consider Schwartz distributions, i.e. continuous functionals on the space ${\cal{S}}(\R^n)$
of rapidly decreasing infinitely differentiable functions on $\R^n$. 
The Fourier transform of a distribution $f$ is defined by $\langle\hat{f}, \phi\rangle= \langle f, \hat{\phi} \rangle$ for
every test function $\phi \in {\cal{S}}(\R^n).$ For any even distribution $f$, we have $(\hat{f})^\wedge
= (2\pi)^n f$.

If $K$ is a convex body  and $0<p<n,$
then $\|\cdot\|_K^{-p}$  is a locally integrable function on $\R^n$ and represents a distribution acting by integration. 
Suppose that $K$ is infinitely smooth, i.e. $\|\cdot\|_K\in C^\infty(S^{n-1})$ is an infinitely differentiable 
function on the sphere. Then by \cite[Lemma 3.16]{K1}, the Fourier transform of $\|\cdot\|_K^{-p}$  
is an extension of some function $g\in C^\infty(S^{n-1})$ to a homogeneous function of degree
$-n+p$ on $\R^n.$ When we write $\left(\|\cdot\|_K^{-p}\right)^\wedge(\xi),$ we mean $g(\xi),\ \xi \in S^{n-1}.$

For $f\in C^\infty(S^{n-1})$ and $0<p<n$, we denote by 
$$(f\cdot r^{-p})(x) = f(x/|x|_2) |x|_2^{-p}$$
the extension of $f$ to a homogeneous function of degree $-p$ on $\R^n.$
Again by  \cite[Lemma 3.16]{K1}, there exists $g\in C^\infty(S^{n-1})$ such that
$$(f\cdot r^{-p})^\wedge = g\cdot r^{-n+p}.$$

If $K,L$ are infinitely smooth convex bodies, the following spherical version of Parseval's
formula was proved in \cite{K4} (see also \cite[Lemma 3.22]{K1}):  for any $p\in (-n,0)$
\begin{equation}\label{parseval}
\int_{S^{n-1}} \left(\|\cdot\|_K^{-p}\right)^\wedge(\xi) \left(\|\cdot\|_L^{-n+p}\right)^\wedge(\xi) =
(2\pi)^n \int_{S^{n-1}} \|x\|_K^{-p} \|x\|_L^{-n+p}\ dx.
\end{equation}

It was proved in \cite[Theorem 1]{K2} that an origin-symmetric convex body $K$ in $\R^n$ is an
intersection body if and only if the function $\|\cdot\|_K^{-1}$ represents a positive definite 
distribution. If $K$ is infinitely smooth, this means that the function $(\|\cdot\|_K^{-1})^\wedge$
is non-negative on the sphere.

We also need a result from \cite{K7} (see also \cite[Theorem 3.8]{K1}) 
expressing volume of central hyperplane sections
in terms of the Fourier transform. For any origin-symmetric star
body $K$ in $\R^n,$ the distribution $(\|\cdot\|_K^{-n+1})^\wedge$ is a continuous function
on the sphere extended to a homogeneous function of degree -1 on the whole of $\R^n,$
and for every $\xi\in S^{n-1},$
\begin{equation} \label{sect-Fourier}
|K\cap \xi^\bot| = \frac 1{\pi(n-1)} (\|\cdot\|_K^{-n+1})^\wedge(\xi).
\end{equation}
In particular, if $K=B_2^n$ and $|\cdot|_2$ is the Euclidean norm, then for every $\xi\in S^{n-1}$
\begin{equation}\label{eucl1}
(|\cdot|_2^{-n+1})^\wedge(\xi) = \pi(n-1)|B_2^{n-1}|.
\end{equation}

\begin{lemma} \label{remain1} 
Let $K$ be an origin-symmetric infinitely smooth convex body in $\R^n.$ Then
$$ \int_{S^{n-1}} \left(\|\cdot\|_K^{-1}\right)^\wedge(\xi) d\xi
\le \frac{(2\pi)^n n}{\pi(n-1)} c_{n,1} |K|^{1/n},$$
 \end{lemma}
\pf 
By (\ref{eucl1}),
Parseval's formula, H\"older's  inequality, polar formula for volume (\ref{polar}) and $|S^{n-1}|=n|B_2^n|$, we get
$$\int_{S^{n-1}} \left(\|\cdot\|_K^{-1}\right)^\wedge(\xi)d\xi$$
$$= \frac{1}{\pi(n-1)\left|B_2^{n-1}\right|} 
\int_{S^{n-1}}  \left(\|\cdot\|_K^{-1}\right)^\wedge(\xi)  \left(|\cdot|_2^{-n+1}\right)^\wedge(\xi) $$
$$=\frac{(2\pi)^n}{\pi(n-1)\left| B_2^{n-1} \right| } \int_{S^{n-1}} \|\theta\|_K^{-1}\ d\theta$$ 
$$ \le  \frac{(2\pi)^n}{\pi(n-1)\left|B_2^{n-1}\right|} \left|S^{n-1}\right|^{\frac{n-1}n} \left(\int_{S^{n-1}} \|\theta\|_K^{-n}\ d\theta\right)^{\frac1n}$$
$$
=  \frac{(2\pi)^n}{\pi(n-1)\left|B_2^{n-1}\right|} \left|S^{n-1}\right|^{\frac{n-1}n} n^{1/n}|K|^{1/n}= 
\frac{(2\pi)^n n}{\pi(n-1)} c_{n,1} |K|^{1/n}. \qed
$$
\bigbreak
The following theorem provides examples where the minimal measure of sections
behaves in a different way from the case of volume.
Note that every non-intersection body can be approximated in
the radial metric by infinitely smooth non-intersection bodies with strictly positive curvature; see 
\cite[Lemma 4.10]{K1}. Different examples of convex bodies that are not intersection bodies
(in dimensions five and higher, as in dimensions up to four such examples do not exist) can
be found in \cite[Chapter 4]{K1}. In particular, the unit balls of the spaces $\ell_q^n,\ q>2,\ n\ge 5$
are not intersection bodies.

\begin{theorem}\label{sepsectcounter} Suppose that $L$ is an infinitely smooth origin-symmetric
convex body in $\R^n$ with strictly positive curvature that is not an intersection body. Then for small enough $\e>0$ 
there exists an origin-symmetric convex body $K$ in $\R^n,\ K\subset L,$ 
such that
$$|K\cap \xi^\bot| \le |L\cap \xi^\bot|-\e,\qquad \forall \xi\in S^{n-1},$$
but
$$|K|^{\frac {n-1}n} > |L|^{\frac {n-1}n} - c_{n,1}\ \e.$$
Note that $c_{n,1}\in (\frac 1{\sqrt{e}},1).$
\end{theorem}

\pf Since $L$ is infinitely smooth, the Fourier transform of $\|\cdot\|_L^{-1}$
is a continuous function on the sphere $S^{n-1}.$ Also, $L$ is not an intersection body,
so $\left(\|\cdot\|_L^{-1}\right)^\wedge < 0$ on an open set $\Omega\subset S^{n-1}.$
Let $\phi\in C^\infty(S^{n-1})$ be an even non-negative, not identically zero, infinitely smooth function
on $S^{n-1}$ with support in $\Omega\cup -\Omega.$ Extend $\phi$ to an even homogeneous
 of degree -1 function $\phi\cdot r^{-1}$ on $\R^n\setminus \{0\}.$ The Fourier transform of this
function in the sense of distributions is $\psi\cdot r^{-n+1}$ where $\psi$ is an infinitely smooth 
function on the sphere.

Let $\e$ be a number such that $|B_2^{n-1}| \|\theta\|_L^{-n+1}> \e>0$ for every $\theta\in S^{n-1}.$
Define a star body $K$ by
\begin{equation}\label{newbody}
\|\theta\|_K^{-n+1}= \|\theta\|_L^{-n+1} - \delta\psi(\theta) - \frac{\e}{|B_2^{n-1}|},\qquad \forall \theta\in S^{n-1},
\end{equation}
where $\delta>0$ is small enough so that for every $\theta$
$$|\delta\psi(\theta)|<\min\left\{\|\theta\|_L^{-n+1}-\frac{\e}{|B_2^{n-1}|},\ \frac{\e}{|B_2^{n-1}|}\right\}.$$ 
The latter condition implies that $K\subset L.$ Since $L$ has strictly positive curvature, by an argument
from \cite[p. 96]{K1}, we can make $\e, \delta$ smaller (if necessary) to ensure that the body $K$ is convex.

Now we extend the functions in (\ref{newbody}) from the sphere to $\R^n\setminus \{0\}$
as homogeneous functions of degree $-n+1$ and apply the Fourier transform. We get
that for every $\xi\in S^{n-1}$
\begin{equation}\label{four0}
\left(\|\cdot\|_K^{-n+1}\right)^{\wedge}(\xi) = \left(\|\cdot\|_L^{-n+1}\right)^\wedge(\xi) - (2\pi)^n \delta \phi(\xi)
- \pi(n-1)\e.
\end{equation}
Here, we used (\ref{eucl1}) to compute the last term.
By (\ref{four0}), (\ref{sect-Fourier}) and the fact that the function $\phi$ is non-negative,
\begin{equation}\label{sect1}
|K\cap \xi^\bot| = |L\cap \xi^\bot| - \frac{(2\pi)^n}{\pi(n-1)} \delta \phi(\xi) - \e \le |L\cap \xi^\bot|-\e.
\end{equation}
Multiplying both sides of (\ref{four0}) by $(\left\|\cdot\|_L^{-1}\right)^\wedge(\xi),$
integrating over $S^{n-1}$ and using Parseval's formula on the sphere,  we get
$$(2\pi)^n \int_{S^{n-1}} \|\theta\|_L^{-1}\ \|\theta\|_K^{-n+1} d\theta$$
$$ =(2\pi)^n n|L| - (2\pi)^n \delta \int_{S^{n-1}} \phi(\theta) (\left\|\cdot\|_L^{-1}\right)^\wedge(\theta) d\theta$$
$$-\pi(n-1) \e \int_{S^{n-1}} (\left\|\cdot\|_L^{-1}\right)^\wedge(\theta) d\theta.$$
Since $\phi$ is a non-negative function supported in $\Omega,$ where $(\left\|\cdot\|_L^{-1}\right)^\wedge$ is negative,
the latter equality implies
$$(2\pi)^n n|L| - \pi(n-1) \e \int_{S^{n-1}} (\left\|\cdot\|_L^{-1}\right)^\wedge(\theta) d\theta $$
$$
< (2\pi)^n \int_{S^{n-1}} \|\theta\|_L^{-1}\ \|\theta\|_K^{-n+1} d\theta 
$$
$$\le (2\pi)^n \left(\int_{S^{n-1}} \|\theta\|_K^{-n} d\theta \right)^{\frac{n-1}n}
\left(\int_{S^{n-1}}\|\theta\|_L^{-n}d\theta\right)^{\frac 1n}$$$$= (2\pi)^n n |L|^{\frac 1n}|K|^{\frac {n-1}{n}}.$$

Combining the latter inequality with the estimate of Lemma \ref{remain1}, we get the result.\qed
\bigbreak

\begin{co}
Suppose that $L$ is an infinitely smooth origin-symmetric
convex body in $\R^n$ with strictly positive curvature that is not an intersection body. 
Then there exists an even continuous function 
$g\ge 0$ on $L$ so that
\begin{equation} \label{meas-counter}
\int_L g < \frac n{n-1} c_{n,1} |L|^{1/n}\min_{\xi\in S^{n-1}} \int_{L\cap \xi^\bot} g .
\end{equation}
\end{co}
\pf By Theorem \ref{sepsectcounter} there exist $\e>0$ and an origin-symmetric convex body $K\subset L$
such that
$$\e=\min_{\xi\in S^{n-1}} \left(|L\cap \xi^\bot| -  |K\cap \xi^\bot|\right),$$
but
$$|L|^{\frac{n-1}n} - |K|^{\frac{n-1}n} < c_{n,1}\e.$$
Combining these and applying the Mean Value Theorem to the function $t\to t^{\frac {n-1}n}$
$$c_{n,1} \min_{\xi\in S^{n-1}} \left(|L\cap \xi^\bot| -  |K\cap \xi^\bot|\right) > |L|^{\frac{n-1}n} - |K|^{\frac{n-1}n}$$
$$ \ge \frac {n-1}n |L|^{-1/n} \left(|L| - |K|\right).$$
The latter shows that $g_0=\chi_{L\setminus K},$ the indicator function of the set $L\setminus K,$
satisfies (\ref{meas-counter}). By simple approximation one can get (\ref{meas-counter})
with a continuous function $g.$\qed

\bigbreak
{\bf Acknowledgement.} I wish to thank the US National Science Foundation for support through 
grant DMS-1265155.

\end{document}